\documentclass[11pt]{article}
\usepackage{a4wide}

\newcommand{\ba}{\begin{array}}
\newcommand{\ea}{\end{array}}
\newcommand{\be}{\begin{equation}}
\newcommand{\ee}{\end{equation}}
\newcommand{\la}{\label}
\newcommand{\bea}{\begin{eqnarray}}
\newcommand{\eea}{\end{eqnarray}}
\renewcommand{\l}{\left}
\renewcommand{\r}{\right}
\newcommand{\n}{\nonumber}
\newcommand{\nn}{\nonumber \\}
\newcommand{\ds}{\displaystyle}
\newcommand{\ndots}{n=0,1,2,\ldots}
\newcommand{\ch}{\choose}
\renewcommand{\a}{\alpha}
\renewcommand{\b}{\beta}
\newcommand{\G}{\Gamma}
\renewcommand{\L}{L_n^{(\a)}(x)}
\newcommand{\KL}{L_n^{\a,M}(x)}
\newcommand{\GL}{L_n^{\a,M,N}(x)}
\renewcommand{\P}{P_n^{(\a,\b)}(x)}
\newcommand{\GP}{P_n^{\a,\b,M,N}(x)}
\newcommand{\SP}{P_n^{(\a,\a)}(x)}
\newcommand{\SGP}{P_n^{\a,\a,M,M}(x)}
\newcommand{\set}[1]{\left\{#1\right\}_{n=0}^{\infty}}
\newcommand{\hyp}[5]{\mbox{}_{#1}F_{#2}
\left(\left.{#3 \atop #4}\right|#5\right)}

\begin{document}

\title{Inversion methods for finding differential equations for generalized
Jacobi polynomials}
\author{J.~Koekoek \and R.~Koekoek}
\date{}
\maketitle

\begin{abstract}
We look for differential equations satisfied by the generalized Jacobi
polynomials $\set{\GP}$ which are orthogonal on the interval $[-1,1]$ with
respect to the weight function
$$\frac{\G(\a+\b+2)}{2^{\a+\b+1}\G(\a+1)\G(\b+1)}(1-x)^{\a}(1+x)^{\b}+
M\delta(x+1)+N\delta(x-1),$$
where $\a>-1$, $\b>-1$, $M\ge 0$ and $N\ge 0$.

In order to find explicit formulas for the coefficients of these
differential equations we have to solve systems of equations of the form
$$\sum_{i=1}^{\infty}A_i(x)D^i\P=F_n(x),\;n=1,2,3,\ldots,$$
where the coefficients $\l\{A_i(x)\r\}_{i=1}^{\infty}$ are independent of
$n$. This system of equations has a unique solution given by
$$A_i(x)=2^i\sum_{j=1}^i\frac{\a+\b+2j+1}{(\a+\b+j+1)_{i+1}}
P_{i-j}^{(-\a-i-1,-\b-i-1)}(x)F_j(x),\;i=1,2,3,\ldots.$$
This is a consequence of the Jacobi inversion formula
\bea & &\sum_{k=j}^i\frac{\a+\b+2k+1}{(\a+\b+k+j+1)_{i-j+1}}\times{}\nn
& &{}\hspace{1cm}{}\times
P_{i-k}^{(-\a-i-1,-\b-i-1)}(x)P_{k-j}^{(\a+j,\b+j)}(x)=\delta_{ij},
\;j\le i,\;i,j=0,1,2,\ldots,\n\eea
which is proved in this report.
\end{abstract}

\vfill

\begin{tabular}{ll}
Keywords : & Differential equations, Generalized Jacobi polynomials,\\
& Inversion formulas.\\[5mm]
\multicolumn{2}{l}{1991 Mathematics Subject Classification :
Primary 33C45 ; Secondary 34A35.}
\end{tabular}

\newpage

\section{Introduction}

In \cite{Koorn} T.H.~Koornwinder introduced the generalized Jacobi
polynomials $\set{\GP}$ which are orthogonal on the interval $[-1,1]$ with
respect to the weight function
$$\frac{\G(\a+\b+2)}{2^{\a+\b+1}\G(\a+1)\G(\b+1)}(1-x)^{\a}(1+x)^{\b}+
M\delta(x+1)+N\delta(x-1),$$
where $\a>-1$, $\b>-1$, $M\ge 0$ and $N\ge 0$. As a limit case he also found
the generalized Laguerre polynomials $\set{\KL}$ which are orthogonal on the
interval $[0,\infty)$ with respect to the weight function
$$\frac{1}{\G(\a+1)}x^{\a}e^{-x}+M\delta(x),$$
where $\a>-1$ and $M\ge 0$. These generalized Jacobi polynomials and
generalized Laguerre polynomials are related by the limit
$$\KL=\lim_{\b\rightarrow\infty}P_n^{\a,\b,0,M}\l(1-\frac{2x}{\b}\r).$$

In \cite{DV} we proved that for $M>0$ the generalized Laguerre polynomials
satisfy a unique differential equation of the form
$$M\sum_{i=0}^{\infty}a_i(x)y^{(i)}(x)+xy''(x)+(\a+1-x)y'(x)+ny(x)=0,$$
where $\l\{a_i(x)\r\}_{i=0}^{\infty}$ are continuous functions on the real
line and $\l\{a_i(x)\r\}_{i=1}^{\infty}$ are independent of the degree $n$.
In \cite{Bav} H.~Bavinck found a new method to obtain the main result of
\cite{DV}. This inversion method was found in a similar way as was done in
\cite{Char} in the case of generalizations of the Charlier polynomials. See
also section~5 for more details. In \cite{Soblag} we used this inversion
method to find all differential equations of the form
\bea\la{DVLag}& &M\sum_{i=0}^{\infty}a_i(x)y^{(i)}(x)+
N\sum_{i=0}^{\infty}b_i(x)y^{(i)}(x)+{}\nn
& &\hspace{1cm}{}+MN\sum_{i=0}^{\infty}c_i(x)y^{(i)}(x)+
xy''(x)+(\a+1-x)y'(x)+ny(x)=0,\eea
where the coefficients $\l\{a_i(x)\r\}_{i=1}^{\infty}$,
$\l\{b_i(x)\r\}_{i=1}^{\infty}$ and $\l\{c_i(x)\r\}_{i=1}^{\infty}$ are
independent of $n$ and the coefficients $a_0(x)$, $b_0(x)$ and $c_0(x)$ are
independent of $x$, satisfied by the Sobolev-type Laguerre polynomials
$\set{\GL}$ which are orthogonal with respect to the inner product
$$<f,g>\;=\frac{1}{\G(\a+1)}\int_0^{\infty}x^{\a}e^{-x}f(x)g(x)dx+Mf(0)g(0)+
Nf'(0)g'(0),$$
where $\a>-1$, $M\ge 0$ and $N\ge 0$. These Sobolev-type Laguerre
polynomials $\set{\GL}$ are generalizations of the generalized Laguerre
polynomials $\set{\KL}$. In fact we have
$$L_n^{\a,M,0}(x)=\KL\;\mbox{ and }\;L_n^{\a,0}(x)=\L.$$

In this report we will prove an inversion formula involving the classical
Jacobi polynomials which can be used to find differential equations of the
form
\bea\la{DVJac}& &M\sum_{i=0}^{\infty}a_i(x)y^{(i)}(x)+
N\sum_{i=0}^{\infty}b_i(x)y^{(i)}(x)+
MN\sum_{i=0}^{\infty}c_i(x)y^{(i)}(x)+{}\nn
& &{}\hspace{1cm}
{}+(1-x^2)y''(x)+\l[\b-\a-(\a+\b+2)x\r]y'(x)+n(n+\a+\b+1)y(x)=0,\eea
where the coefficients $\l\{a_i(x)\r\}_{i=1}^{\infty}$,
$\l\{b_i(x)\r\}_{i=1}^{\infty}$ and $\l\{c_i(x)\r\}_{i=1}^{\infty}$ are
independent of $n$ and the coefficients $a_0(x)$, $b_0(x)$ and $c_0(x)$ are
independent of $x$, satisfied by the generalized Jacobi polynomials
$\set{\GP}$. As an example we apply the special case $\b=\a$ of this
inversion formula to solve the systems of equations obtained in
\cite{Symjac}.

The inversion formula for the Charlier polynomials obtained in \cite{Char}
(see also section~5) was also used in \cite{SobChar} to find difference
operators with Sobolev-type Charlier polynomials as eigenfunctions.

In \cite{Meixner} H.~Bavinck and H.~van~Haeringen used similar inversion
formulas to find difference equations for generalized Meixner polynomials.

\section{The classical Laguerre and Jacobi polynomials}

In this section we list the definitions and some properties of the classical
Laguerre and Jacobi polynomials which we will use in this report. For details
the reader is referred to \cite{Chihara}, \cite{AS} and \cite{Szego}.

The classical Laguerre polynomials $\set{\L}$ can be defined by
\be\la{defLag}\L=\sum_{k=0}^n(-1)^k{n+\a \ch n-k}\frac{x^k}{k!},\;\ndots\ee
for all $\a$. Their generating function is given by
\be\la{genLag}(1-t)^{-\a-1}\exp\l(\frac{xt}{t-1}\r)=
\sum_{n=0}^{\infty}\L t^n\ee
and for all $n\in\{0,1,2,\ldots\}$ we have
\be\la{diffLag}D^i\L=(-1)^iL_{n-i}^{(\a+i)}(x),\;i=0,1,2,\ldots,n,\ee
where $D=\ds\frac{d}{dx}$ denotes the differentiation operator. The Laguerre
polynomials satisfy the linear second order differential equation
\be\la{dvLag}xy''(x)+(\a+1-x)y'(x)+ny(x)=0.\ee
It is well-known that
\be\la{formLag}\frac{x^n}{n!}=
\sum_{k=0}^n(-1)^k{n+\a \ch n-k}L_k^{(\a)}(x),\;\ndots.\ee
This formula can easily be proved by using definition (\ref{defLag}) and
changing the order of summation as follows
\bea\sum_{k=0}^n(-1)^k{n+\a \ch n-k}L_k^{(\a)}(x)&=&\sum_{k=0}^n(-1)^k
{n+\a \ch n-k}\sum_{j=0}^k(-1)^j{k+\a \ch k-j}\frac{x^j}{j!}\nn
&=&\sum_{j=0}^n\sum_{k=0}^{n-j}(-1)^k{n+\a \ch n-j-k}{j+k+\a \ch k}
\frac{x^j}{j!}\nn
&=&\sum_{j=0}^n{n+\a \ch n-j}\frac{x^j}{j!}\sum_{k=0}^{n-j}(-1)^k{n-j \ch k}
=\frac{x^n}{n!},\;\ndots.\n\eea

The classical Jacobi polynomials $\set{\P}$ can be defined by
\bea\la{defJac1}\P&=&\sum_{k=0}^n\frac{(n+\a+\b+1)_k}{k!}
\frac{(\a+k+1)_{n-k}}{(n-k)!}\l(\frac{x-1}{2}\r)^k,\;\ndots\\
&=&\la{defJac2}(-1)^n\sum_{k=0}^n\frac{(-n-k-\a-\b)_k}{k!}
\frac{(-n-\a)_{n-k}}{(n-k)!}\l(\frac{x-1}{2}\r)^k,\;\ndots\\
&=&\la{defJac3}2^{-n}\sum_{k=0}^n{n+\a \ch n-k}{n+\b \ch k}
(x-1)^k(x+1)^{n-k},\;\ndots\eea
for all $\a$ and $\b$. For all $n\in\{0,1,2,\ldots\}$ we have
\be\la{diffJac}D^i\P=\frac{(n+\a+\b+1)_i}{2^i}P_{n-i}^{(\a+i,\b+i)}(x),
\;i=0,1,2,\ldots,n.\ee
The Jacobi polynomials satisfy the symmetry formula
\be\la{symJac}P_n^{(\a,\b)}(-x)=(-1)^nP_n^{(\b,\a)}(x),\;\ndots\ee
and the linear second order differential equation
\be\la{dvJac}(1-x^2)y''(x)+\l[\b-\a-(\a+\b+2)x\r]y'(x)
+n(n+\a+\b+1)y(x)=0.\ee
Further we have for $\a+\b+1>0$ (compare with \cite{Luke}, page 277, formula
(30))
\be\la{formJac*}\l(\frac{1-x}{2}\r)^n=\sum_{k=0}^n
\frac{(-n)_k(\a+k+1)_{n-k}(\a+\b+2k+1)}{(\a+\b+k+1)_{n+1}}P_k^{(\a,\b)}(x),
\;\ndots.\ee
This formula is much less known than formula (\ref{formLag}) for the
Laguerre polynomials. However, the proof is quite similar. In section~6 we
will prove a much more general formula.

We remark that (\ref{formJac*}) can be written in a more general form as
\bea\la{formJac} & &\sum_{k=0}^n
\frac{(-n)_k(\a+\b+1)_k(\a+k+1)_{n-k}(\a+\b+2k+1)}{\G(\a+\b+n+k+2)}
P_k^{(\a,\b)}(x)\nn
&=&\frac{1}{\G(\a+\b+1)}\l(\frac{1-x}{2}\r)^n,\;\ndots,\eea
which is valid for all $\a$ and $\b$.

The Jacobi polynomials $\set{\P}$ and the Laguerre polynomials $\set{\L}$
are related by the limit
\be\la{limit}\L=\lim_{\b\rightarrow\infty}
P_n^{(\a,\b)}\l(1-\frac{2x}{\b}\r).\ee
We remark that if we replace $x$ by $\ds 1-\frac{2x}{\b}$ in (\ref{formJac*}),
multiply by $\b^n$ and let $\b$ tend to infinity in the complex plane along
the halfline where $\a+\b$ is real and $\a+\b+1>0$ we obtain (\ref{formLag})
by using (\ref{limit}) and the fact that we have for all
$n\in\{0,1,2,\ldots\}$
$$(-n)_k(\a+k+1)_{n-k}=(-1)^k{n+\a \ch n-k}n!,k=0,1,2,\ldots,n.$$

Finally we derive some formulas involving Jacobi polynomials which we will
need in this report. By using definition (\ref{defJac1}) we easily obtain
\bea & &n\P-(n+\a)P_{n-1}^{(\a,\b+1)}(x)\nn
&=&n\sum_{k=0}^n\frac{(n+\a+\b+1)_k}{k!}\frac{(\a+k+1)_{n-k}}{(n-k)!}
\l(\frac{x-1}{2}\r)^k+{}\nn
& &{}\hspace{1cm}{}-(n+\a)\sum_{k=0}^{n-1}\frac{(n+\a+\b+1)_k}{k!}
\frac{(\a+k+1)_{n-k-1}}{(n-k-1)!}\l(\frac{x-1}{2}\r)^k\nn
&=&\sum_{k=0}^n\frac{(n+\a+\b+1)_k}{k!}\frac{(\a+k+1)_{n-k}}{(n-k)!}
\l[n-(n-k)\r]\l(\frac{x-1}{2}\r)^k\nn
&=&\sum_{k=0}^n\frac{(n+\a+\b+1)_k}{k!}\frac{(\a+k+1)_{n-k}}{(n-k)!}
(x-1)D\l(\frac{x-1}{2}\r)^k,\;n=1,2,3,\ldots.\n\eea
Hence
\be\la{relJac1}(x-1)D\P=n\P-(n+\a)P_{n-1}^{(\a,\b+1)}(x),\;n=1,2,3,\ldots.\ee
By using the symmetry property (\ref{symJac}) we also have
\be\la{relJac2}(x+1)D\P=n\P+(n+\b)P_{n-1}^{(\a+1,\b)}(x),\;n=1,2,3,\ldots.\ee
Another easy consequence of definition (\ref{defJac1}) is
\be\la{relJac3}P_n^{(\a+1,\b)}(x)-P_n^{(\a,\b+1)}(x)=
P_{n-1}^{(\a+1,\b+1)}(x),\;n=1,2,3,\ldots.\ee

If $\b=\a$ we can combine the relations (\ref{relJac1}), (\ref{relJac2}) and
(\ref{relJac3}) to find that
\be\la{relsymJac1}2xD\SP=2n\SP+(n+\a)P_{n-2}^{(\a+1,\a+1)}(x),
\;n=2,3,4,\ldots.\ee
This formula turns out to be very useful in section~7. In that section we
will also need the following well-known relation for ultraspherical
polynomials (see for instance \cite{Szego})
\bea\la{relsymJac2} & &(n+2\a+1)(n+2\a+2)P_n^{(\a+1,\a+1)}(x)
-(n+\a)(n+\a+1)P_{n-2}^{(\a+1,\a+1)}(x)\nn
& &{}\hspace{1cm}{}=2(n+\a+1)(2n+2\a+1)\SP,\;n=2,3,4,\ldots.\eea
This formula can be proved by straightforward calculations as follows. As
before we find by using definition (\ref{defJac1})
\bea & &(n+2\a+1)(n+2\a+2)P_n^{(\a+1,\a+1)}(x)
-(n+\a)(n+\a+1)P_{n-2}^{(\a+1,\a+1)}(x)\nn
&=&(n+2\a+1)(n+2\a+2)\sum_{k=0}^n\frac{(n+2\a+3)_k}{k!}\frac{(\a+k+2)_{n-k}}{(n-k)!}
\l(\frac{x-1}{2}\r)^k+{}\nn
& &{}\hspace{1cm}{}-(n+\a)(n+\a+1)\sum_{k=0}^{n-2}\frac{(n+2\a+1)_k}{k!}
\frac{(\a+k+2)_{n-k-2}}{(n-k-2)!}\l(\frac{x-1}{2}\r)^k\nn
&=&\sum_{k=0}^n\frac{(n+2\a+1)_k}{k!}\frac{(\a+k+2)_{n-k}}{(n-k)!}\times{}\nn
& &{}\hspace{1cm}{}\times
\l[(n+2\a+k+1)(n+2\a+k+2)-(n-k)(n-k-1)\r]\l(\frac{x-1}{2}\r)^k\nn
&=&2(2n+2\a+1)\sum_{k=0}^n\frac{(n+2\a+1)_k}{k!}
\frac{(\a+k+1)_{n-k+1}}{(n-k)!}\l(\frac{x-1}{2}\r)^k\nn
&=&2(n+\a+1)(2n+2\a+1)\SP,\;n=2,3,4,\ldots.\n\eea
Another useful relation involving ultraspherical polynomials is
\be\la{relsymJac3}(\a+1)P_n^{(\a+1,\a+1)}(x)-(n+\a+1)\SP=
\frac{1}{4}(n+\a+1)(1-x^2)P_{n-2}^{(\a+2,\a+2)}(x),\ee
which holds for $n=2,3,4,\ldots$. This can easily be shown by using
definition (\ref{defJac3}) as follows
\bea & &(\a+1)P_n^{(\a+1,\a+1)}(x)-(n+\a+1)\SP\nn
&=&2^{-n}\sum_{k=0}^n\l[(\a+1){n+\a+1 \ch n-k}{n+\a+1 \ch k}+{}\r.\nn
& &{}\hspace{2cm}\l.{}-(n+\a+1){n+\a \ch n-k}{n+\a \ch k}\r]
(x-1)^k(x+1)^{n-k}\nn
&=&-(n+\a+1)2^{-n}\sum_{k=1}^{n-1}{n+\a \ch n-k-1}{n+\a \ch k-1}
(x-1)^k(x+1)^{n-k}\nn
&=&(n+\a+1)(1-x^2)2^{-n}\sum_{k=0}^{n-2}{n+\a \ch n-k-2}{n+\a \ch k}
(x-1)^k(x+1)^{n-k-2}\nn
&=&\frac{1}{4}(n+\a+1)(1-x^2)P_{n-2}^{(\a+2,\a+2)}(x),\;n=2,3,4,\ldots.\n\eea

\section{Some summation formulas}

In this section we will derive some summation formulas which will be used in
this report. We start with a partial sum of a ${}_1F_0$ hypergeometric series
at the point $1$. In fact we have
\be\la{binsum}\sum_{k=0}^n{a+k \ch k}=\sum_{k=0}^n\frac{(a+1)_k}{k!}=
\frac{(a+2)_n}{n!}={n+a+1 \ch n},\;\ndots.\ee
This formula appeared as a special case of lemma~4 in \cite{Soblag}. Further
we have the well-known Vandermonde summation formula
\be\la{Van}\hyp{2}{1}{-n,b}{c}{1}=\frac{(c-b)_n}{(c)_n},\;
(c)_n\ne 0,\;\ndots.\ee
By using (\ref{Van}) we find that
\bea\la{nul} & &\sum_{k=0}^n\frac{(-n)_k(b)_k}{(b+n+1)_kk!}(b+2k)\nn
&=&b\,\hyp{2}{1}{-n,b+1}{b+n+1}{1}-
\frac{nb}{b+n+1}\,\hyp{2}{1}{-n+1,b+1}{b+n+2}{1}\nn
&=&b\,\frac{(n)_n}{(b+n+1)_n}-
\frac{nb}{b+n+1}\frac{(n+1)_{n-1}}{(b+n+2)_{n-1}}=0,
\;(b+n+1)_n\ne 0,\;n=1,2,3,\ldots.\hspace{5mm}{}\eea
We remark that this can be written in a more general form as
\be\la{nulalg}\sum_{k=0}^n\frac{(-n)_k(b)_k}{\G(b+n+k+1)k!}(b+2k)=0,
\;n=1,2,3,\ldots,\ee
which is valid for all $b$. This formula will be used in section~6.

We also have the well-known Saalsch\"utz summation formula
\be\la{Saal}\hyp{3}{2}{-n,a,n+b+c-1}{b,a+c}{1}=
\frac{(b-a)_n(c)_n}{(b)_n(a+c)_n},\;(b)_n(a+c)_n\ne 0,\;\ndots.\ee
Note that this implies that
\be\la{breuk}\frac{(c)_n}{(b)_n}=
\hyp{3}{2}{-n,\frac{1}{2}(b-c),n+b+c-1}{b,\frac{1}{2}(b+c)}{1},
\;(b)_n((b+c)/2)_n\ne 0,\;\ndots.\ee
This can be used to show that for
$(\frac{1}{2}b+\frac{1}{2})_n(b+n+1)_n(b-c+1)_n\ne 0$
\be\la{sum}\sum_{k=0}^n\frac{(-n)_k(b)_k(c)_k}{(b+n+1)_k(b-c+1)_kk!}(b+2k)
=\frac{(b)_{n+1}(\frac{1}{2}b-c+\frac{1}{2})_n}
{(b-c+1)_n(\frac{1}{2}b+\frac{1}{2})_n},\;\ndots.\ee
So we suppose that $n\in\{0,1,2,\ldots\}$ and that
$(\frac{1}{2}b+\frac{1}{2})_n(b+n+1)_n(b-c+1)_n\ne 0$. Then we use
(\ref{breuk}) and change the order of summation to obtain
\bea & &\sum_{k=0}^n\frac{(-n)_k(b)_k(c)_k}{(b+n+1)_k(b-c+1)_kk!}(b+2k)\nn
&=&\sum_{k=0}^n\frac{(-n)_k(b)_k}{(b+n+1)_kk!}(b+2k)\,
\hyp{3}{2}{-k,\frac{1}{2}b-c+\frac{1}{2},b+k}
{b-c+1,\frac{1}{2}b+\frac{1}{2}}{1}\nn
&=&\sum_{k=0}^n\frac{(-n)_k(b)_k}{(b+n+1)_kk!}(b+2k)\,
\sum_{m=0}^k\frac{(-k)_m(\frac{1}{2}b-c+\frac{1}{2})_m(b+k)_m}
{(b-c+1)_m(\frac{1}{2}b+\frac{1}{2})_mm!}\nn
&=&\sum_{m=0}^n\sum_{k=0}^{n-m}
\frac{(-n)_{m+k}(-m-k)_m(b)_{2m+k}(\frac{1}{2}b-c+\frac{1}{2})_m}
{(b+n+1)_{m+k}(b-c+1)_m(\frac{1}{2}b+\frac{1}{2})_m(m+k)!\,m!}(b+2m+2k)\nn
&=&\sum_{m=0}^n(-1)^m\frac{(-n)_m(b)_{2m}(\frac{1}{2}b-c+\frac{1}{2})_m}
{(b+n+1)_m(b-c+1)_m(\frac{1}{2}b+\frac{1}{2})_mm!}\times{}\nn
& &{}\hspace{2cm}{}\times\sum_{k=0}^{n-m}\frac{(-n+m)_k(b+2m)_k}
{(b+n+m+1)_kk!}(b+2m+2k).\n\eea
Hence, by using (\ref{nul}) we conclude that
\bea & &\sum_{k=0}^n\frac{(-n)_k(b)_k(c)_k}{(b+n+1)_k(b-c+1)_kk!}(b+2k)\nn
&=&(-1)^n\frac{(-n)_n(b)_{2n}(\frac{1}{2}b-c+\frac{1}{2})_n}
{(b+n+1)_n(b-c+1)_n(\frac{1}{2}b+\frac{1}{2})_nn!}(b+2n)\nn
&=&\frac{(b)_{2n+1}(\frac{1}{2}b-c+\frac{1}{2})_n}
{(b+n+1)_n(b-c+1)_n(\frac{1}{2}b+\frac{1}{2})_n}
=\frac{(b)_{n+1}(\frac{1}{2}b-c+\frac{1}{2})_n}
{(b-c+1)_n(\frac{1}{2}b+\frac{1}{2})_n},\;\ndots,\n\eea
which proves (\ref{sum}).

\section{The systems of equations}

Let $\a>-1$. The Sobolev-type Laguerre polynomials $\set{\GL}$ can be
written as
$$\GL=A_0\L+A_1D\L+A_2D^2\L,\;\ndots,$$
where the coefficients $A_0$, $A_1$ and $A_2$ are given by
$$\l\{\ba{l}\ds A_0=1+M{n+\a \ch n-1}+\frac{n(\a+2)-(\a+1)}{(\a+1)(\a+3)}
N{n+\a \ch n-2}+{}\\
\ds\hspace{5cm}{}+\frac{MN}{(\a+1)(\a+2)}{n+\a \ch n-1}{n+\a+1 \ch n-2}\\  \\
\ds A_1=M{n+\a \ch n}+\frac{n-1}{\a+1}N{n+\a \ch n-1}+
\frac{2MN}{(\a+1)^2}{n+\a \ch n}{n+\a+1 \ch n-2}\\  \\
\ds A_2=\frac{N}{\a+1}{n+\a \ch n-1}+\frac{MN}{(\a+1)^2}{n+\a \ch n}
{n+\a+1 \ch n-1}.\ea\r.$$
For details concerning these Sobolev-type Laguerre polynomials and their
definition the reader is referred to \cite{Thesis} and \cite{SIAM}. Since
the classical Laguerre polynomials $\set{\L}$ satisfy the differential
equation (\ref{dvLag}) it is quite reasonable to look for differential
equations of the form (\ref{DVLag}) for these Sobolev-type Laguerre
polynomials $\set{\GL}$ in view of this definition and the fact that
$L_n^{\a,0,0}(x)=\L$. In \cite{Soblag} it is shown that this leads to eight
systems of equations for the coefficients $\l\{a_i(x)\r\}_{i=0}^{\infty}$,
$\l\{b_i(x)\r\}_{i=0}^{\infty}$ and $\l\{c_i(x)\r\}_{i=0}^{\infty}$. In
order to find these coefficients we have to solve systems of equations which
are of the form
$$\sum_{i=1}^{\infty}A_i(x)D^{i+k}\L=F_n(x),\;n=k+1,k+2,k+3,\ldots,$$
where $k\in\{0,1,2,\ldots\}$ and the coefficients
$\l\{A_i(x)\r\}_{i=1}^{\infty}$ are independent of $n$. In \cite{Soblag} it
is pointed out that this system of equations has a unique solution given by
$$A_i(x)=(-1)^{i+k}\sum_{j=1}^iL_{i-j}^{(-\a-i-k-1)}(-x)F_{j+k}(x),
\;i=1,2,3,\ldots.$$
This is an easy consequence of the Laguerre inversion formula
\be\la{invLag}\sum_{k=j}^iL_{i-k}^{(-\a-i-1)}(-x)L_{k-j}^{(\a+j)}(x)=
\delta_{ij},\;j\le i,\;i,j=0,1,2,\ldots,\ee
which was found by H.~Bavinck in \cite{Bav}. For more details the reader is
referred to \cite{Bav} and \cite{Soblag}. See also section~5 of this report.

Now we take $\a>-1$ and $\b>-1$. The generalized Jacobi polynomials
$\set{\GP}$ can be written as
$$\GP=A_0\P+\l[A_1(1-x)-A_2(1+x)\r]D\P,\;\ndots,$$
where the coefficients $A_0$, $A_1$ and $A_2$ are given by
$$\l\{\ba{l}\ds A_0=1+
M\frac{\ds {n+\b \ch n-1}{n+\a+\b+1 \ch n}}{\ds {n+\a \ch n}}+
N\frac{\ds {n+\a \ch n-1}{n+\a+\b+1 \ch n}}{\ds {n+\b \ch n}}+{}\\  \\
\ds\hspace{7cm}{}+
MN\frac{(\a+\b+2)^2}{(\a+1)(\b+1)}{n+\a+\b+1 \ch n-1}^2\\ \\
\ds A_1=\frac{M}{\a+\b+1}
\frac{\ds {n+\b \ch n}{n+\a+\b \ch n}}{\ds {n+\a \ch n}}
+\frac{MN}{\a+1}{n+\a+\b \ch n-1}{n+\a+\b+1 \ch n}\\ \\
\ds A_2=\frac{N}{\a+\b+1}
\frac{\ds {n+\a \ch n}{n+\a+\b \ch n}}{\ds {n+\b \ch n}}
+\frac{MN}{\b+1}{n+\a+\b \ch n-1}{n+\a+\b+1 \ch n}.\ea\r.$$
Here we used the same definition as in \cite{Koorn}, but in a slightly
different notation. The case $\a+\b+1=0$ must be understood by continuity.
In view of this definition and the fact that the classical Jacobi
polynomials $\set{\P}$ satisfy the differential equation (\ref{dvJac}) it
is quite natural to look for differential equations of the form
(\ref{DVJac}) satisfied by these generalized Jacobi polynomials as was
already pointed out in \cite{JCAM}. Again this leads to eight systems of
equations for the coefficients $\l\{a_i(x)\r\}_{i=0}^{\infty}$,
$\l\{b_i(x)\r\}_{i=0}^{\infty}$ and $\l\{c_i(x)\r\}_{i=0}^{\infty}$. In
order to find these coefficients we have to solve systems of equations which
are of the form
$$\sum_{i=1}^{\infty}A_i(x)D^i\P=F_n(x),\;n=1,2,3,\ldots,$$
where the coefficients $\l\{A_i(x)\r\}_{i=1}^{\infty}$ are independent of
$n$. This system of equations has a unique solution given by
$$A_i(x)=2^i\sum_{j=1}^i\frac{\a+\b+2j+1}{(\a+\b+j+1)_{i+1}}
P_{i-j}^{(-\a-i-1,-\b-i-1)}(x)F_j(x),\;i=1,2,3,\ldots.$$
This is a consequence of the Jacobi inversion formula
\bea\la{invJac}& &\sum_{k=j}^i
\frac{\a+\b+2k+1}{(\a+\b+k+j+1)_{i-j+1}}\times{}\nn
& &{}\hspace{1cm}{}\times
P_{i-k}^{(-\a-i-1,-\b-i-1)}(x)P_{k-j}^{(\a+j,\b+j)}(x)=\delta_{ij},
\;j\le i,\;i,j=0,1,2,\ldots,\eea
which will be proved in this report. Again, the case $\a+\b+1=0$ must be
understood by continuity. We remark that if we replace $x$ by
$\ds 1-\frac{2x}{\b}$ in (\ref{invJac}), multiply by $\b^{i-j}$ and let $\b$
tend to infinity along the positive real axis we obtain the Laguerre
inversion formula (\ref{invLag}) by using (\ref{limit}).

In \cite{Symjac} we found all differential equations of the form
\be\la{DVSGP}M\sum_{i=0}^{\infty}a_i(x)y^{(i)}(x)+
(1-x^2)y''(x)-2(\a+1)xy'(x)+n(n+2\a+1)y(x)=0,\ee
where $\l\{a_i(x)\r\}_{i=0}^{\infty}$ are continuous functions on the real
line and $\l\{a_i(x)\r\}_{i=1}^{\infty}$ are independent of $n$, satisfied
by the symmetric generalized ultraspherical polynomials $\set{\SGP}$ defined
by
$$\SGP=C_0\SP-C_1xD\SP,\;\ndots,$$
where
$$\l\{\ba{l}\ds C_0=1+\frac{2Mn}{\a+1}{n+2\a+1 \ch n}
+4M^2{n+2\a+1 \ch n-1}^2\\
\\
\ds C_1=\frac{2M}{2\a+1}{n+2\a \ch n}
+\frac{2M^2}{\a+1}{n+2\a \ch n-1}{n+2\a+1 \ch n}.\ea\r.$$
We remark that these polynomials form a special case ($\b=\a$ and $N=M$) of
the generalized Jacobi polynomials $\set{\GP}$, but the differential
equation (\ref{DVSGP}) has a very special form without a $M^2$-part. This
differential equation will appear not to be a special case of the
differential equation of the form (\ref{DVJac}) for the generalized Jacobi
polynomials, since the $MN$-part will not vanish if we take $\b=\a$ and
$N=M$. We aim to give a proof of this in a future publication. In section~7
of this report we will apply the special case $\b=\a$ of the Jacobi inversion
formula (\ref{invJac}) to solve the systems of equations obtained in
\cite{Symjac} as an example.

\section{The inversion formulas}

In \cite{Char} H.~Bavinck and R.~Koekoek found the following inversion
formula involving Charlier polynomials
\be\la{invChar}\sum_{k=j}^iC_{i-k}^{(-a)}(-x)C_{k-j}^{(a)}(x)=
\delta_{ij},\;j\le i,\;i,j=0,1,2,\ldots.\ee
This formula is an easy consequence of the generating function (see for
instance \cite{AS})
$$e^{-at}(1+t)^x=\sum_{n=0}^{\infty}C_n^{(a)}(x)t^n.$$
In fact we have
\bea 1&=&e^{-at}(1+t)^xe^{at}(1+t)^{-x}=
\sum_{k=0}^{\infty}C_k^{(a)}(x)t^k\sum_{m=0}^{\infty}C_m^{(-a)}(-x)t^m\nn
&=&\sum_{n=0}^{\infty}\l(\sum_{k=0}^nC_k^{(a)}(x)C_{n-k}^{(-a)}(-x)\r)t^n.
\n\eea
Hence
$$\sum_{k=0}^nC_k^{(a)}(x)C_{n-k}^{(-a)}(-x)=\l\{\ba{ll}1, & n=0\\ \\
0, & n=1,2,3,\ldots.\ea\r.$$
Now (\ref{invChar}) easily follows by taking $n=i-j$ and shifting the
summation index. This formula was also used in \cite{SobChar} to find
difference operators with Sobolev-type Charlier polynomials as
eigenfunctions. In \cite{Meixner} a similar formula involving Meixner
polynomials was used to find difference equations for generalized Meixner
polynomials.

Formula (\ref{invChar}) can be interpreted as follows. If we define the
matrix $T=(t_{ij})_{i,j=0}^n$ with entries
$$t_{ij}=\l\{\ba{ll}C_{i-j}^{(a)}(x), & j\le i\\ \\0, & j>i,\ea\r.$$
then this matrix $T$ is a triangular matrix with determinant $1$ and the
inverse $U$ of this matrix is given by $T^{-1}=U=(u_{ij})_{i,j=0}^n$ with
entries
$$u_{ij}=\l\{\ba{ll}C_{i-j}^{(-a)}(-x), & j\le i\\ \\0, & j>i.\ea\r.$$
Therefore we call (\ref{invChar}) an inversion formula.

In the same way we find by using the generating function (\ref{genLag})
for the Laguerre polynomials
$$\sum_{k=j}^iL_{i-k}^{(\a)}(x)L_{k-j}^{(-\a-2)}(-x)=\delta_{ij},\;j\le i,
\;i,j=0,1,2,\ldots.$$
However, this formula cannot be used to solve systems of equations of the
form
$$\sum_{i=1}^{\infty}A_i(x)D^i\L=F_n(x),\;n=1,2,3,\ldots$$
in view of the parametershift in (\ref{diffLag}).

In \cite{Bav} H.~Bavinck used a slightly different method to find the
Laguerre inversion formula (\ref{invLag}) from the generating function
(\ref{genLag}) for the Laguerre polynomials. In fact we have
\bea (1-t)^{i-j-1}&=&(1-t)^{-\a-j-1}\exp\l(\frac{xt}{t-1}\r)
(1-t)^{\a+i}\exp\l(\frac{-xt}{t-1}\r)\nn
&=&\sum_{k=0}^{\infty}L_k^{(\a+j)}(x)t^k
\sum_{m=0}^{\infty}L_m^{(-\a-i-1)}(-x)t^m\nn
&=&\sum_{n=0}^{\infty}
\l(\sum_{k=0}^nL_k^{(\a+j)}(x)L_{n-k}^{(-\a-i-1)}(-x)\r)t^n.\n\eea
This implies, by comparing the coefficients of $t^{i-j}$ on both sides, that
$$\sum_{k=0}^{i-j}L_k^{(\a+j)}(x)L_{i-j-k}^{(-\a-i-1)}(-x)=\delta_{ij},
\;j\le i,\;i,j=0,1,2,\ldots,$$
which is equivalent to (\ref{invLag}).

Formula (\ref{invLag}) can be interpreted as follows. If we define the
matrix $T=(t_{ij})_{i,j=0}^n$ with entries
$$t_{ij}=\l\{\ba{ll}L_{i-j}^{(\a+j)}(x), & j\le i\\ \\0, & j>i,\ea\r.$$
then this matrix $T$ is a triangular matrix with determinant $1$ and the
inverse $U$ of this matrix is given by $T^{-1}=U=(u_{ij})_{i,j=0}^n$ with
entries
$$u_{ij}=\l\{\ba{ll}L_{i-j}^{(-\a-i-1)}(-x), & j\le i\\ \\0, & j>i.\ea\r.$$

It is possible to generalize the Laguerre inversion formula to
\be\la{geninvLag}\sum_{k=0}^nL_k^{(\a+p_n)}(x)L_{n-k}^{(-\a-q_n)}(-x)=
\frac{(p_n-q_n+2)_n}{n!},\;\ndots,\ee
where $p_n$ and $q_n$ are arbitrary and even may depend on $n$. In order to
have an inversion formula we have to choose $p_n$ and $q_n$ such that
$$(p_n-q_n+2)_n=0,\;n=1,2,3,\ldots,$$
hence
$$p_n-q_n\in\{-n-1,-n,\ldots,-3,-2\},\;n=1,2,3,\ldots.$$
Note that the endpoint-cases $p_n-q_n=-2$ and $p_n-q_n=-n-1$ correspond to
the earlier mentioned inversion formulas.

To prove (\ref{geninvLag}) we use (\ref{genLag}) to obtain
$$\L=\frac{1}{n!}\l.D_t^n\l[(1-t)^{-\a-1}
\exp\l(\frac{xt}{t-1}\r)\r]\r|_{t=0},\;\ndots,$$
where $\ds D_t=\frac{d}{dt}$ denotes differentiation with respect to $t$.
Hence by using Leibniz' rule we find
\bea & &\sum_{k=0}^nL_k^{(\a+p_n)}(x)L_{n-k}^{(-\a-q_n)}(-x)\nn
&=&\sum_{k=0}^n\frac{1}{k!}\l.D_t^k\l[(1-t)^{-\a-p_n-1}
\exp\l(\frac{xt}{t-1}\r)\r]\r|_{t=0}\times{}\nn
& &{}\hspace{2cm}{}\times\frac{1}{(n-k)!}
\l.D_t^{n-k}\l[(1-t)^{\a+q_n-1}\exp\l(\frac{-xt}{t-1}\r)\r]\r|_{t=0}\nn
&=&\frac{1}{n!}\l.D_t^n\l[(1-t)^{q_n-p_n-2}\r]\r|_{t=0}=
\frac{(p_n-q_n+2)_n}{n!},\;\ndots,\n\eea
which proves (\ref{geninvLag}).

In case of the Jacobi polynomials the above methods seem not to be
applicable. However, in the next section we will give a proof of the Jacobi
inversion formula (\ref{invJac}).

\section{Proof of the Jacobi inversion formula}

In this section we will prove that
\bea\la{algJac}& &\sum_{k=0}^n
\frac{(\a+\b+2k+1)(\a+\b+1)_k}{\G(\a+\b+n+k+2)}P_k^{(\a,\b)}(x)
P_{n-k}^{(-n-\a-1,-n-\b-1)}(y)\nn
&=&\frac{1}{\G(\a+\b+1)}\frac{1}{n!}\l(\frac{x-y}{2}\r)^n,
\;n=0,1,2,\ldots,\eea
which holds for all $\a$ and $\b$.

Note that (\ref{formJac}) is a special case of (\ref{algJac}) since
\bea P_{n-k}^{(-n-\a-1,-n-\b-1)}(1)&=&\frac{(-n-\a)_{n-k}}{(n-k)!}=
(-1)^{n-k}\frac{(\a+k+1)_{n-k}}{(n-k)!}\nn
&=&\frac{(-1)^n}{n!}(-n)_k(\a+k+1)_{n-k},\;k=0,1,2,\ldots,n\n\eea
for all $n\in\{0,1,2,\ldots\}$.

By taking $y=x$ in (\ref{algJac}) we easily obtain
\bea\la{inv}& &\sum_{k=0}^n
\frac{(\a+\b+2k+1)(\a+\b+1)_k}{\G(\a+\b+n+k+2)}\times{}\nn
& &{}\hspace{1cm}{}\times P_k^{(\a,\b)}(x)P_{n-k}^{(-n-\a-1,-n-\b-1)}(x)=
\l\{\ba{ll}\ds\frac{1}{\G(\a+\b+1)}, & n=0\\ \\0, & n=1,2,3,\ldots\ea\r.\eea
for all $\a$ and $\b$. If we take $n=i-j$ in (\ref{inv}) and shift the
summation index we find
\bea & &\sum_{k=j}^i
\frac{(\a+\b+2k-2j+1)(\a+\b+1)_{k-j}}{\G(\a+\b+i-2j+k+2)}\times{}\nn
& &{}\hspace{1cm}{}\times
P_{i-k}^{(-i+j-\a-1,-i+j-\b-1)}(x)P_{k-j}^{(\a,\b)}(x)=
\frac{\delta_{ij}}{\G(\a+\b+1)},\;j\le i,\;i,j=0,1,2,\ldots.\n\eea
For $\a$ and $\b$ real with $\a+\b+1>-1$ we now obtain (\ref{invJac}) by
shifting both $\a$ and $\b$ by $j$.

Note that (\ref{algJac}) for $y=-x$ in a similar way leads to
\bea\la{spec}& &\sum_{k=j}^i
\frac{\a+\b+2k+1}{(\a+\b+k+j+1)_{i-j+1}}\times{}\nn
& &{}\hspace{1cm}{}\times
P_{i-k}^{(-\a-i-1,-\b-i-1)}(-x)P_{k-j}^{(\a+j,\b+j)}(x)=
\frac{x^{i-j}}{(i-j)!},\;j\le i,\;i,j=0,1,2,\ldots,\eea
which will also be used in the next section.

In order to prove (\ref{algJac}) we start with the left-hand side, apply
definition (\ref{defJac1}) to $P_k^{(\a,\b)}(x)$ and definition
(\ref{defJac2}) to $P_{n-k}^{(-n-\a-1,-n-\b-1)}(y)$ and change the order of
summation to obtain
\bea & &\sum_{k=0}^n
\frac{(\a+\b+2k+1)(\a+\b+1)_k}{\G(\a+\b+n+k+2)}P_k^{(\a,\b)}(x)
P_{n-k}^{(-n-\a-1,-n-\b-1)}(y)\nn
&=&\sum_{k=0}^n\sum_{i=0}^k\sum_{j=0}^{n-k}(-1)^{n-k}
\frac{(\a+\b+2k+1)(\a+\b+1)_k}{\G(\a+\b+n+k+2)}\frac{(\a+\b+k+1)_i}{i!}
\frac{(\a+i+1)_{k-i}}{(k-i)!}\nn
& &{}\hspace{1cm}{}\times
\frac{(\a+\b+n+k-j+2)_j}{j!}\frac{(\a+k+1)_{n-k-j}}{(n-k-j)!}
\l(\frac{x-1}{2}\r)^i\l(\frac{y-1}{2}\r)^j\nn
&=&\sum_{i=0}^n\sum_{k=i}^n\sum_{j=0}^{n-k}(-1)^{n-k}\times{}\nn
& &{}\hspace{1cm}{}\times
\frac{(\a+\b+2k+1)(\a+\b+1)_{i+k}(\a+i+1)_{n-i-j}}
{\G(\a+\b+n+k-j+2)\,i!\,(k-i)!\,j!\,(n-k-j)!}
\l(\frac{x-1}{2}\r)^i\l(\frac{y-1}{2}\r)^j\nn
&=&\sum_{i=0}^n\sum_{k=0}^{n-i}\sum_{j=0}^{n-i-k}(-1)^{n-i-k}\times{}\nn
& &{}\hspace{1cm}{}\times
\frac{(\a+\b+2i+2k+1)(\a+\b+1)_{2i+k}(\a+i+1)_{n-i-j}}
{\G(\a+\b+n+i+k-j+2)\,i!\,k!\,j!\,(n-i-k-j)!}
\l(\frac{x-1}{2}\r)^i\l(\frac{y-1}{2}\r)^j\nn
&=&\sum_{i=0}^n\sum_{j=0}^{n-i}\sum_{k=0}^{n-i-j}(-1)^{n-i}
\l(\frac{x-1}{2}\r)^i\l(\frac{y-1}{2}\r)^j\times{}\nn
& &{}\hspace{1cm}{}\times
\frac{(\a+\b+2i+2k+1)(\a+\b+1)_{2i+k}(\a+i+1)_{n-i-j}(-n+i+j)_k}
{\G(\a+\b+n+i-j+k+2)\,i!\,j!\,k!\,(n-i-j)!}\nn
&=&\sum_{i=0}^n\sum_{j=0}^{n-i}(-1)^{n-i}
\frac{(\a+\b+1)_{2i}(\a+i+1)_{n-i-j}}{i!\,j!\,(n-i-j)!}
\l(\frac{x-1}{2}\r)^i\l(\frac{y-1}{2}\r)^j\times{}\nn
& &{}\hspace{1cm}{}\times\sum_{k=0}^{n-i-j}
\frac{(-n+i+j)_k(\a+\b+2i+1)_k}{\G(\a+\b+n+i-j+k+2)k!}(\a+\b+2i+2k+1),
\;\ndots.\n\eea
Now we use (\ref{nulalg}) to obtain
\bea & &\sum_{k=0}^n
\frac{(\a+\b+2k+1)(\a+\b+1)_k}{\G(\a+\b+n+k+2)}P_k^{(\a,\b)}(x)
P_{n-k}^{(-n-\a-1,-n-\b-1)}(y)\nn
&=&\sum_{i=0}^n(-1)^{n-i}\frac{(\a+\b+1)_{2i}}{i!\,(n-i)!}
\l(\frac{x-1}{2}\r)^i\l(\frac{y-1}{2}\r)^{n-i}\frac{\a+\b+2i+1}{\G(\a+\b+2i+2)}\nn
&=&\frac{1}{\G(\a+\b+1)}\frac{1}{n!}\sum_{i=0}^n{n \ch i}
\l(\frac{x-1}{2}\r)^i\l(\frac{1-y}{2}\r)^{n-i}\nn
&=&\frac{1}{\G(\a+\b+1)}\frac{1}{n!}\l(\frac{x-y}{2}\r)^n,\;\ndots,\n\eea
which proves (\ref{algJac}).

\section{The symmetric generalized ultraspherical polynomials}

Let $\a>-1$. In \cite{Symjac} we found the coefficients
$\l\{a_i(x)\r\}_{i=0}^{\infty}$ of the differential equation (\ref{DVSGP})
for the symmetric generalized ultraspherical polynomials $\set{\SGP}$. In
order to do this we had to solve the following two systems of equations for
the coefficients $\l\{a_i(x)\r\}_{i=0}^{\infty}$~:
\be\la{sys1}\sum_{i=0}^{\infty}a_i(x)D^i\SP=
\frac{4}{2\a+1}{n+2\a \ch n}D^2\SP\ee
and
\be\la{sys2}\sum_{i=0}^{\infty}ia_i(x)D^i\SP+
x\sum_{i=0}^{\infty}a_i(x)D^{i+1}\SP=4{n+2\a+1 \ch n-1}D^2\SP\ee
for $\ndots$, where the coefficients $\l\{a_i(x)\r\}_{i=0}^{\infty}$ are
continuous functions on the real line and $\l\{a_i(x)\r\}_{i=1}^{\infty}$
are independent of $n$. Now we suppose that $a_0(x):=a_0(n,\a)$ is
independent of $x$ as we did in \cite{Soblag}. Then it is clear (see for
instance lemma 1 in \cite{Soblag}) that $a_i(x)$ must be a polynomial in $x$
of degree at most $i$ for each $i=1,2,3,\ldots$. In \cite{Symjac} we showed
that the solution for $\l\{a_i(x)\r\}_{i=0}^{\infty}$ is not unique. In fact
it was shown that
\be\la{anul}a_0(x):=a_0(n,\a)=a_0(1,\a)b_0(n,\a)+c_0(n,\a),\;\ndots\ee
and that
\be\la{a}a_i(x)=a_0(1,\a)b_i(x)+c_i(x),\;i=1,2,3,\ldots,\ee
where $a_0(1,\a)$ is arbitrary and
\be\la{bnul}b_0(n,\a)=\frac{1}{2}\l[1-(-1)^n\r],\;\ndots,\ee
\be\la{cnul}c_0(n,\a)=4(2\a+3){n+2\a+2 \ch n-2},\;\ndots,\ee
\be\la{b}b_i(x)=\frac{2^{i-1}}{i!}(-x)^i,\;i=1,2,3,\ldots,\ee
\be\la{c}c_1(x)=0\;\mbox{ and }\;c_i(x)=(2\a+3)(1-x^2)
\frac{2^i}{i!}P_{i-2}^{(\a-i+3,\a-i+3)}(x),\;i=2,3,4,\ldots.\ee

In this section we will give an alternative proof of this by using the
Jacobi inversion formula (\ref{invJac}). Note that the special case $\b=\a$
of the Jacobi inversion formula (\ref{invJac}) reads
\bea\la{invsym}& &\sum_{k=j}^i
\frac{2\a+2k+1}{(2\a+k+j+1)_{i-j+1}}\times{}\nn
& &{}\hspace{1cm}{}\times
P_{i-k}^{(-\a-i-1,-\a-i-1)}(x)P_{k-j}^{(\a+j,\a+j)}(x)=\delta_{ij},
\;j\le i,\;i,j=0,1,2,\ldots.\eea
If we apply this inversion formula to the system of equations
\be\la{sysalg}\sum_{i=1}^{\infty}A_i(x)D^i\SP=F_n(x),\;n=1,2,3,\ldots,\ee
then we find
\be\la{solalg}A_i(x)=2^i\sum_{j=1}^i\frac{2\a+2j+1}{(2\a+j+1)_{i+1}}
P_{i-j}^{(-\a-i-1,-\a-i-1)}(x)F_j(x),\;i=1,2,3,\ldots.\ee

The special case $\b=\a$ of (\ref{spec}) reads
\bea\la{specsym}& &\sum_{k=j}^i
\frac{2\a+2k+1}{(2\a+k+j+1)_{i-j+1}}\times{}\nn
& &{}\hspace{1cm}{}\times
P_{i-k}^{(-\a-i-1,-\a-i-1)}(-x)P_{k-j}^{(\a+j,\a+j)}(x)=
\frac{x^{i-j}}{(i-j)!},\;j\le i,\;i,j=0,1,2,\ldots.\eea

By considering (\ref{sys1}) and (\ref{sys2}) for $n=0$ and $n=1$ we conclude
that $a_0(0,\a)=0$, $a_0(1,\a)$ is arbitrary and $a_1(x)=-a_0(1,\a)x$. For
$n=2,3,4,\ldots$ it turns out to be more convenient to use another system of
equations instead of (\ref{sys2}). By using (\ref{relsymJac1}) we find for
$i=0,1,2,\ldots$
\bea & &iD^i\SP+xD^{i+1}\SP=D^i\l[xD\SP\r]\nn
&=&nD^i\SP+\frac{1}{2}(n+\a)D^iP_{n-2}^{(\a+1,\a+1)}(x),
\;n=2,3,4,\ldots.\n\eea
Combining (\ref{sys1}) and (\ref{sys2}) we now obtain
\bea\sum_{i=0}^{\infty}a_i(x)D^iP_{n-2}^{(\a+1,\a+1)}(x)&=&
\frac{2}{n+\a}\l[4{n+2\a+1 \ch n-1}-\frac{4n}{2\a+1}{n+2\a \ch n}\r]D^2\SP\nn
&=&\frac{8}{n+\a}\l[{n+2\a+1 \ch n-1}-{n+2\a \ch n-1}\r]D^2\SP\nn
&=&\frac{8}{n+\a}{n+2\a \ch n-2}D^2\SP,\;n=2,3,4,\ldots.\n\eea
So we conclude that (\ref{sys2}) for $n=2,3,4,\ldots$ may be replaced by
\be\la{sys3}\sum_{i=0}^{\infty}a_i(x)D^iP_{n-2}^{(\a+1,\a+1)}(x)=
\frac{8}{n+\a}{n+2\a \ch n-2}D^2\SP,\;n=2,3,4,\ldots.\ee
Note that for $n=2$ this implies that $a_0(2,\a)=4(2\a+3)$.

Since $a_i(x)$ must be a polynomial in $x$ of degree at most $i$ for each
$i=1,2,3,\ldots$ we may write
$$a_i(x)=k_ix^i+\mbox{ lower order terms },\;i=1,2,3,\ldots.$$
By comparing the coefficients of highest degree in (\ref{sys1}) and
(\ref{sys3}) we find by using (\ref{defJac1})~:
$$\frac{a_0(n,\a)}{n!}+\sum_{i=1}^n\frac{k_i}{(n-i)!}=0,\;n=1,2,3,\ldots$$
and
$$\frac{a_0(n,\a)}{(n-2)!}+\sum_{i=1}^{n-2}\frac{k_i}{(n-i-2)!}=
4(2n+2\a-1){n+2\a \ch n-2}\frac{1}{(n-2)!},\;n=3,4,5,\ldots.$$
Since $k_i$ is independent of $n$ for $i=1,2,3,\ldots$ and
$a_0(2,\a)=4(2\a+3)$ we conclude that
\be\la{verschil}a_0(n+2,\a)-a_0(n,\a)=4(2n+2\a+3){n+2\a+2 \ch n},
\;\ndots,\ee
where $a_0(0,\a)=0$ and $a_0(1,\a)$ is arbitrary. Hence
$$a_0(2n,\a)-a_0(0,\a)=4\sum_{k=0}^{n-1}{2k+2\a+2 \ch 2k}(4k+2\a+3),
\;n=1,2,3,\ldots$$
and
$$a_0(2n+1,\a)-a_0(1,\a)=4\sum_{k=0}^{n-1}{2k+2\a+3 \ch 2k+1}(4k+2\a+5),
\;n=1,2,3,\ldots.$$
Since we have by using (\ref{binsum})
\bea & &\sum_{k=0}^{n-1}{2k+2\a+2 \ch 2k}(4k+2\a+3)=
\sum_{k=0}^{n-1}\frac{(2\a+3)_{2k+1}}{(2k)!}
+\sum_{k=0}^{n-1}\frac{2k(2\a+3)_{2k}}{(2k)!}\nn
&=&\sum_{k=0}^{2n-2}\frac{(2\a+3)_{k+1}}{k!}=
(2\a+3)\frac{(2\a+5)_{2n-2}}{(2n-2)!}=(2\a+3){2n+2\a+2 \ch 2n-2},
\;n=1,2,3,\ldots\n\eea
and
\bea & &\sum_{k=0}^{n-1}{2k+2\a+3 \ch 2k+1}(4k+2\a+5)=
\sum_{k=0}^{n-1}\frac{(2\a+3)_{2k+2}}{(2k+1)!}
+\sum_{k=0}^{n-1}\frac{(2\a+3)_{2k+1}}{(2k)!}\nn
&=&\sum_{k=0}^{2n-1}\frac{(2\a+3)_{k+1}}{k!}=
(2\a+3)\frac{(2\a+5)_{2n-1}}{(2n-1)!}=(2\a+3){2n+2\a+3 \ch 2n-1},
\;n=1,2,3,\ldots,\n\eea
we now conclude that (\ref{anul}), (\ref{bnul}) and (\ref{cnul}) hold.

The systems of equations (\ref{sys1}) and (\ref{sys3}) lead to
\be\la{sys4}\sum_{i=1}^{\infty}a_i(x)D^i\SP=
\frac{4}{2\a+1}{n+2\a \ch n}D^2\SP-a_0(n,\a)\SP\ee
for $\ndots$ and
\be\la{sys5}\sum_{i=1}^{\infty}a_i(x)D^iP_{n-2}^{(\a+1,\a+1)}(x)=
\frac{8}{n+\a}{n+2\a \ch n-2}D^2\SP-a_0(n,\a)P_{n-2}^{(\a+1,\a+1)}(x)\ee
for $n=2,3,4,\ldots$.

First we remark that (\ref{sys4}) is true for $n=0$ and $n=1$ since
$a_0(0,\a)=0$ and $a_1(x)=-a_0(1,\a)x$. Then we will show that every
solution of (\ref{sys5}) also satisfies (\ref{sys4}). Suppose that
$\l\{a_i(x)\r\}_{i=1}^{\infty}$ is a solution of (\ref{sys5}). Now we use
(\ref{diffJac}), (\ref{relsymJac2}), (\ref{verschil}) and the fact that
$\l\{a_i(x)\r\}_{i=1}^{\infty}$ are independent of $n$ to obtain for
$n=2,3,4,\ldots$
\bea & &2(n+\a+1)(2n+2\a+1)\sum_{i=1}^{\infty}a_i(x)D^i\SP\nn
&=&(n+2\a+1)(n+2\a+2)\sum_{i=1}^{\infty}a_i(x)D^iP_n^{(\a+1,\a+1)}(x)+{}\nn
& &{}\hspace{1cm}{}
-(n+\a)(n+\a+1)\sum_{i=1}^{\infty}a_i(x)D^iP_{n-2}^{(\a+1,\a+1)}(x)\nn
&=&(n+2\a+1)(n+2\a+2)\times{}\nn
& &{}\hspace{2cm}{}\times
\l[\frac{8}{n+\a+2}{n+2\a+2 \ch n}D^2P_{n+2}^{(\a,\a)}(x)
-a_0(n+2,\a)P_n^{(\a+1,\a+1)}(x)\r]+{}\nn
& &{}\hspace{1cm}{}-(n+\a)(n+\a+1)\l[\frac{8}{n+\a}{n+2\a \ch n-2}D^2\SP
-a_0(n,\a)P_{n-2}^{(\a+1,\a+1)}(x)\r]\nn
&=&\frac{2(n+2\a+1)(n+2\a+2)}{n+\a+2}{n+2\a+2 \ch n}\times{}\nn
& &{}\hspace{2cm}{}\times
\l[(n+2\a+3)(n+2\a+4)P_n^{(\a+2,\a+2)}(x)+{}\r.\nn
& &{}\hspace{6cm}{}\l.-2(n+\a+2)(2n+2\a+3)P_n^{(\a+1,\a+1)}(x)\r]+{}\nn
& &{}\hspace{1cm}{}-8(n+\a+1){n+2\a \ch n-2}D^2\SP+{}\nn
& &{}\hspace{1cm}{}-a_0(n,\a)
\l[(n+2\a+1)(n+2\a+2)P_n^{(\a+1,\a+1)}(x)+{}\r.\nn
& &{}\hspace{6cm}{}\l.-(n+\a)(n+\a+1)P_{n-2}^{(\a+1,\a+1)}(x)\r]\nn
&=&2(n+\a+1)(n+2\a+1)(n+2\a+2){n+2\a+2 \ch n}P_{n-2}^{(\a+2,\a+2)}(x)+{}\nn
& &{}\hspace{1cm}{}-8(n+\a+1){n+2\a \ch n-2}D^2\SP+{}\nn
& &{}\hspace{1cm}{}-2(n+\a+1)(2n+2\a+1)a_0(n,\a)\SP\nn
&=&8(n+\a+1)\l[{n+2\a+2 \ch n}-{n+2\a \ch n-2}\r]D^2\SP+{}\nn
& &{}\hspace{1cm}{}-2(n+\a+1)(2n+2\a+1)a_0(n,\a)\SP\nn
&=&2(n+\a+1)(2n+2\a+1)\l[\frac{4}{2\a+1}{n+2\a \ch n}D^2\SP
-a_0(n,\a)\SP\r].\n\eea
Since $\a>-1$ this proves that every solution of (\ref{sys5}) also satisfies
(\ref{sys4}).

Now we will solve (\ref{sys5}). Shifting $n$ by two we may write, since the
coefficients $\l\{a_i(x)\r\}_{i=1}^{\infty}$ are independent of $n$
\be\la{sys6}\sum_{i=1}^{\infty}a_i(x)D^iP_n^{(\a+1,\a+1)}(x)=F_n(x),
\;\ndots,\ee
where
$$F_n(x)=\frac{8}{n+\a+2}{n+2\a+2 \ch n}D^2P_{n+2}^{(\a,\a)}(x)
-a_0(n+2,\a)P_n^{(\a+1,\a+1)}(x).$$
Since $a_0(2,\a)=4(2\a+3)$ we easily find that $F_0(x)=0$. This implies that
the system of equations (\ref{sys6}) is of the form (\ref{sysalg}). So if we
apply the inversion formula (\ref{invsym}) to the system of equations
(\ref{sys6}) we obtain by using (\ref{solalg})
$$a_i(x)=2^i\sum_{j=1}^i\frac{2\a+2j+3}{(2\a+j+3)_{i+1}}
P_{i-j}^{(-\a-i-2,-\a-i-2)}(x)F_j(x),\;i=1,2,3,\ldots.$$
Hence, by using (\ref{anul}) we now conclude that the coefficients
$\l\{a_i(x)\r\}_{i=1}^{\infty}$ can be written in the form (\ref{a}).
Moreover, we find by using (\ref{bnul}), (\ref{symJac}) and (\ref{specsym})
\bea b_i(x)&=&2^{i-1}\sum_{j=1}^i\frac{2\a+2j+3}{(2\a+j+3)_{i+1}}
P_{i-j}^{(-\a-i-2,-\a-i-2)}(x)P_j^{(\a+1,\a+1)}(x)\l[(-1)^j-1\r]\nn
&=&2^{i-1}\sum_{j=0}^i\frac{2\a+2j+3}{(2\a+j+3)_{i+1}}
P_{i-j}^{(-\a-i-2,-\a-i-2)}(x)P_j^{(\a+1,\a+1)}(-x)+{}\nn
& &{}\hspace{1cm}{}-2^{i-1}\sum_{j=0}^i\frac{2\a+2j+3}{(2\a+j+3)_{i+1}}
P_{i-j}^{(-\a-i-2,-\a-i-2)}(x)P_j^{(\a+1,\a+1)}(x)\nn
&=&2^{i-1}\frac{(-x)^i}{i!},\;i=1,2,3,\ldots,\n\eea
which proves (\ref{b}). And by using (\ref{cnul}) and (\ref{diffJac}) we
obtain
$$c_i(x)=2^i\sum_{j=1}^i\frac{2\a+2j+3}{(2\a+j+3)_{i+1}}
P_{i-j}^{(-\a-i-2,-\a-i-2)}(x)G_j(x),\;i=1,2,3,\ldots,$$
where
$$G_j(x)=\frac{4(2\a+3)}{j+\a+2}{j+2\a+4 \ch j}
\l[(\a+2)P_j^{(\a+2,\a+2)}(x)-(j+\a+2)P_j^{(\a+1,\a+1)}(x)\r].$$
It is clear that $G_1(x)=0$, which implies that $c_1(x)=0$. Note that since
$b_1(x)=-x$ this also implies that $a_1(x)=-a_0(1,\a)x$, which agrees with
what we have found before. Now we use (\ref{relsymJac3}) to find
$$G_j(x)=(2\a+3)(1-x^2){j+2\a+4 \ch j}P_{j-2}^{(\a+3,\a+3)}(x),
\;j=2,3,4,\ldots.$$
Hence, for $i=2,3,4,\ldots$ we have
\bea & &c_i(x)=(2\a+3)(1-x^2)2^i\times{}\nn
& &{}\hspace{2cm}{}\times\sum_{j=2}^i\frac{2\a+2j+3}{(2\a+j+3)_{i+1}}
{j+2\a+4 \ch j}P_{i-j}^{(-\a-i-2,-\a-i-2)}(x)P_{j-2}^{(\a+3,\a+3)}(x).\n\eea
Now it remains to show that
\bea\la{cc} & &\sum_{j=2}^i\frac{2\a+2j+3}{(2\a+j+3)_{i+1}}
{j+2\a+4 \ch j}P_{i-j}^{(-\a-i-2,-\a-i-2)}(x)P_{j-2}^{(\a+3,\a+3)}(x)\nn
& &{}\hspace{1cm}{}
=\frac{1}{i!}P_{i-2}^{(\a-i+3,\a-i+3)}(x),\;i=2,3,4,\ldots.\eea
In order to do this we write for $i=2,3,4,\ldots$
\bea & &\sum_{j=2}^i\frac{2\a+2j+3}{(2\a+j+3)_{i+1}}
{j+2\a+4 \ch j}P_{i-j}^{(-\a-i-2,-\a-i-2)}(x)P_{j-2}^{(\a+3,\a+3)}(x)\nn
& &{}\hspace{1cm}{}=\sum_{k=0}^{i-2}\frac{2\a+2k+7}{(2\a+k+5)_{i+1}}
\frac{(2\a+5)_{k+2}}{(k+2)!}P_{i-k-2}^{(-\a-i-2,-\a-i-2)}(x)
P_k^{(\a+3,\a+3)}(x).\n\eea
Now we apply definition (\ref{defJac1}) to $P_k^{(\a+3,\a+3)}(x)$ and
definition (\ref{defJac2}) to $P_{i-k-2}^{(-\a-i-2,-\a-i-2)}(x)$ and change
the order of summation to obtain for $i=2,3,4,\ldots$
\bea & &\sum_{k=0}^{i-2}\frac{2\a+2k+7}{(2\a+k+5)_{i+1}}
\frac{(2\a+5)_{k+2}}{(k+2)!}P_{i-k-2}^{(-\a-i-2,-\a-i-2)}(x)
P_k^{(\a+3,\a+3)}(x)\nn
&=&\sum_{k=0}^{i-2}\sum_{m=0}^{i-k-2}\sum_{n=0}^k(-1)^{i-k-2}\times{}\nn
& &{}\hspace{1cm}{}\times\frac{(2\a+2k+7)(2\a+5)_{k+n+2}(\a+n+4)_{i-m-n-2}}
{(2\a+k+5)_{i-m+1}(k+2)!\,m!\,(i-k-m-2)!\,n!\,(k-n)!}\l(\frac{x-1}{2}\r)^{m+n}\nn
&=&\sum_{m=0}^{i-2}\sum_{n=0}^{i-m-2}\sum_{k=0}^{i-m-n-2}(-1)^{i-n-k}\times{}\nn
& &{}\hspace{1cm}{}\times\frac{(2\a+2n+2k+7)(2\a+5)_{2n+k+2}(\a+n+4)_{i-m-n-2}}
{(2\a+n+k+5)_{i-m+1}(n+k+2)!\,m!\,(i-m-n-k-2)!\,n!\,k!}\l(\frac{x-1}{2}\r)^{m+n}\nn
&=&\sum_{j=0}^{i-2}\sum_{n=0}^j\sum_{k=0}^{i-j-2}(-1)^{i-n-k}\times{}\nn
& &{}\hspace{1cm}{}\times\frac{(2\a+2n+2k+7)(2\a+5)_{n+k}(\a+n+4)_{i-j-2}}
{(2\a+2n+k+7)_{i-j-1}(n+k+2)!\,(j-n)!\,(i-j-k-2)!\,n!\,k!}\l(\frac{x-1}{2}\r)^j\nn
&=&\sum_{j=0}^{i-2}\sum_{n=0}^j(-1)^{i-n}\frac{(2\a+5)_n(\a+n+4)_{i-j-2}}
{(2\a+2n+7)_{i-j-1}(n+2)!\,(j-n)!\,(i-j-2)!\,n!}\l(\frac{x-1}{2}\r)^j\times{}\nn
& &{}\hspace{1cm}{}\times\sum_{k=0}^{i-j-2}
\frac{(-i+j+2)_k(2\a+n+5)_k(2\a+2n+7)_k}{(2\a+2n+i-j+6)_k(n+3)_kk!}
(2\a+2n+2k+7).\n\eea
Hence, by using (\ref{sum}), (\ref{Van}) and definition (\ref{defJac2}) we
obtain for $i=2,3,4,\ldots$
\bea & &\sum_{k=0}^{i-2}\frac{2\a+2k+7}{(2\a+k+5)_{i+1}}
\frac{(2\a+5)_{k+2}}{(k+2)!}P_{i-k-2}^{(-\a-i-2,-\a-i-2)}(x)
P_k^{(\a+3,\a+3)}(x)\nn
&=&\sum_{j=0}^{i-2}\sum_{n=0}^j(-1)^{i-n}\frac{(2\a+5)_n(\a+n+4)_{i-j-2}}
{(2\a+2n+7)_{i-j-1}(n+2)!\,(j-n)!\,(i-j-2)!\,n!}\l(\frac{x-1}{2}\r)^j\times{}\nn
& &{}\hspace{1cm}{}\times\frac{(2\a+2n+7)_{i-j-1}(-\a-1)_{i-j-2}}
{(n+3)_{i-j-2}(\a+n+4)_{i-j-2}}\nn
&=&\sum_{j=0}^{i-2}\sum_{n=0}^j(-1)^{i-n}\frac{(2\a+5)_n(-\a-1)_{i-j-2}}
{(n+i-j)!\,(j-n)!\,(i-j-2)!\,n!}\l(\frac{x-1}{2}\r)^j\nn
&=&(-1)^i\sum_{j=0}^{i-2}\frac{(-\a-1)_{i-j-2}}{(i-j)!\,(i-j-2)!\,j!}
\l(\frac{x-1}{2}\r)^j\hyp{2}{1}{-j,2\a+5}{i-j+1}{1}\nn
&=&\frac{(-1)^i}{i!}\sum_{j=0}^{i-2}
\frac{(-\a-1)_{i-j-2}(i-j-2\a-4)_j}{(i-j-2)!\,j!}\l(\frac{x-1}{2}\r)^j\nn
&=&\frac{(-1)^{i-2}}{i!}\sum_{j=0}^{i-2}\frac{(i-j-2\a-4)_j}{j!}
\frac{(-\a-1)_{i-j-2}}{(i-j-2)!}\l(\frac{x-1}{2}\r)^j
=\frac{1}{i!}P_{i-2}^{(\a-i+3,\a-i+3)}(x),\n\eea
which proves (\ref{cc}).

\section{Some remarks}

Note that we have from definition (\ref{defLag}) for the Laguerre
polynomials that
\be\la{monLag}L_n^{(-n)}(x)=(-1)^n\frac{x^n}{n!},\;\ndots.\ee
Hence, the polynomial $L_n^{(-n)}(x)$ reduces to a monomial of degree $n$
for all $n\in\{0,1,2,\ldots\}$. Definition (\ref{defJac3}) for the Jacobi
polynomials leads to
$$P_n^{(-n,\b)}(x)={n+\b \ch n}\l(\frac{x-1}{2}\r)^n,\;\ndots,$$
which is also a monomial. However, this monomial might reduce to the zero
polynomial. For instance, $P_n^{(-n,-n)}(x)$ equals the zero polynomial for
all $n\in\{1,2,3,\ldots\}$.

Further we remark that if we replace $x$ by $\ds 1-\frac{2x}{\b}$ and $y$ by
$\ds 1-\frac{2y}{\b}$ in (\ref{algJac}), multiply by $\G(\a+\b+1)\b^n$ and
let $\b$ tend to infinity in an appropriate way we obtain by using
(\ref{limit})
\be\la{algLag}\sum_{k=0}^nL_k^{(\a)}(x)L_{n-k}^{(-n-\a-1)}(-y)=
\frac{(y-x)^n}{n!},\;\ndots.\ee
Note that (\ref{formLag}) is a special case of (\ref{algLag}) since
$$L_{n-k}^{(-n-\a-1)}(0)=(-1)^{n-k}\frac{(\a+k+1)_{n-k}}{(n-k)!}=
(-1)^n(-1)^k{n+\a \ch n-k},\;k=0,1,2,\ldots,n$$
for all $n\in\{0,1,2,\ldots\}$. Moreover, note that (\ref{algLag}) is a
special case of the well-known convolution formula for the classical Laguerre
polynomials
$$\sum_{k=0}^nL_k^{(\a)}(x)L_{n-k}^{(\b)}(y)=L_n^{(\a+\b+1)}(x+y),\;\ndots$$
in view of (\ref{monLag}). By using the technique demonstrated at the end of
section~5 this convolution formula can be proved for all $\a$ and $\b$ which
might even depend on $n$.

We remark that, by using the fact that
$$(b/2)_k(b+2k)=b(b/2+1)_k,\;k=0,1,2,\ldots,$$
the formulas (\ref{nul}) and (\ref{sum}) can also be obtained by using
summation formulas for terminating well-poised hypergeometric series (see
for instance formulas (III.9) and (III.26) in \cite{Slater}).

The computation of $a_0(n,\a)$ from (\ref{verschil}) can also be done in the
following way. Note that we have
$$(2n+2\a+3){n+2\a+2 \ch n}=
(2\a+3)\l[{n+2\a+4 \ch n}-{n+2\a+2 \ch n-2}\r],\;\ndots.$$
Hence, instead of (\ref{binsum}) we can also use the telescoping property of
the sums to find that
$$\sum_{k=0}^{n-1}{2k+2\a+2 \ch 2k}(4k+2\a+3)=(2\a+3){2n+2\a+2 \ch 2n-2},
\;n=1,2,3,\ldots$$
and
$$\sum_{k=0}^{n-1}{2k+2\a+3 \ch 2k+1}(4k+2\a+5)=(2\a+3){2n+2\a+3 \ch 2n-1},
\;n=1,2,3,\ldots.$$

Finally, we remark that if we apply the inversion formula (\ref{invsym}) to
the system of equations (\ref{sys4}) instead of (\ref{sys5}) we find for
$i=1,2,3,\ldots$ that
\be\la{b*}b_i(x)=2^{i-1}\sum_{j=1}^i\frac{2\a+2j+1}{(2\a+j+1)_{i+1}}
P_{i-j}^{(-\a-i-1,-\a-i-1)}(x)P_j^{(\a,\a)}(x)\l[(-1)^j-1\r]\ee
and
\bea\la{c*} c_i(x)&=&2^{i+2}\sum_{j=1}^i\frac{2\a+2j+1}{(2\a+j+1)_{i+1}}
P_{i-j}^{(-\a-i-1,-\a-i-1)}(x)\times{}\nn
& &{}\hspace{1cm}{}\times\l[\frac{1}{2\a+1}{j+2\a \ch j}D^2P_j^{(\a,\a)}(x)
-(2\a+3){j+2\a+2 \ch j-2}P_j^{(\a,\a)}(x)\r].\hspace{5mm}{}\eea
From (\ref{b*}) we easily obtain (\ref{b}) in the same way as before by
using (\ref{symJac}), (\ref{invsym}) and (\ref{specsym}). Further we easily
find from (\ref{c*}) that $c_1(x)=0$, but we were not able to derive
(\ref{c}) for $i=2,3,4,\ldots$ from (\ref{c*}).

\vspace{10mm}

Menelaoslaan 4, 5631 LN Eindhoven, The Netherlands

\vspace{5mm}

Delft University of Technology, Faculty of Technical Mathematics and
Informatics,

P.O. Box 5031, 2600 GA Delft, The Netherlands, e-mail :
koekoek@twi.tudelft.nl


\begin{thebibliography}{99}
\bibitem{Bav} {\sc H. Bavinck :} {\em A direct approach to Koekoek's
differential equation for generalized Laguerre polynomials.} Acta
Mathematica Hungarica {\bf 66}, 1995, 247-253.
\bibitem{SobChar} {\sc H. Bavinck :} {\em A difference operator of
infinite order with Sobolev-type Charlier polynomials as eigenfunctions.}
Indagationes Mathematicae {\bf 7}, 1996, 281-291.
\bibitem{Meixner} {\sc H. Bavinck \& H. van Haeringen :} {\em Difference
equations for generalized Meixner polynomials.} Journal of Mathematical
Analysis and Applications {\bf 184}, 1994, 453-463.
\bibitem{Char} {\sc H. Bavinck \& R. Koekoek :} {\em On a difference
equation for generalizations of Charlier polynomials.} Journal of
Approximation Theory {\bf 81}, 1995, 195-206.
\bibitem{Chihara} {\sc T.S. Chihara :} {\em An introduction to orthogonal
polynomials.} Mathematics and Its Applications {\bf 13}, Gordon and Breach,
New York, 1978.
\bibitem{DV} {\sc J. Koekoek \& R. Koekoek :} {\em On a differential equation
for Koornwinder's generalized Laguerre polynomials.} Proceedings of the American
Mathematical Society {\bf 112}, 1991, 1045-1054.
\bibitem{Soblag} {\sc J. Koekoek, R. Koekoek \& H. Bavinck :} {\em On
differential equations for Sobolev-type Laguerre polynomials.} Transactions
of the American Mathematical Society, to appear.
\bibitem{Thesis} {\sc R. Koekoek :} {\em Generalizations of the classical
Laguerre polynomials and some q-analogues.} Delft University of Technology,
Thesis, 1990.
\bibitem{JCAM} {\sc R. Koekoek :} {\em The search for differential equations
for certain sets of orthogonal polynomials.} Journal of Computational and
Applied Mathematics {\bf 49}, 1993, 111-119.
\bibitem{Symjac} {\sc R. Koekoek :} {\em Differential equations for
symmetric generalized ultraspherical polynomials.} Transactions of the
American Mathematical Society {\bf 345}, 1994, 47-72.
\bibitem{SIAM} {\sc R. Koekoek \& H.G. Meijer :} {\em A generalization
of Laguerre polynomials.} SIAM Journal on Mathematical Analysis {\bf 24},
1993, 768-782.
\bibitem{AS} {\sc R. Koekoek \& R.F. Swarttouw :} {\em The Askey-scheme of
hypergeometric orthogonal polynomials and its $q$-analogue.} Delft
University of Technology, Faculty of Technical Mathematics and Informatics,
report no. {\bf 94-05}, 1994.
\bibitem{Koorn} {\sc T.H. Koornwinder :} {\em Orthogonal polynomials with weight
function $(1-x)^{\a}(1+x)^{\b}+M\delta(x+1)+N\delta(x-1)$.}
Canadian Mathematical Bulletin {\bf 27}(2), 1984, 205-214.
\bibitem{Luke} {\sc Y.L. Luke :} {\em The special functions and their
approximations.} Volume I. Academic Press, San Diego, 1969.
\bibitem{Slater} {\sc L.J. Slater :} {\em Generalized hypergeometric
functions.} Cambridge University Press, Cambridge, 1966.
\bibitem{Szego} {\sc G. Szeg{\"o} :} {\em Orthogonal polynomials.}
American Mathematical Society Colloquium Publications {\bf 23} (1939),
Fourth edition, Providence, Rhode Island, 1975.
\end{thebibliography}
\end{document}